
\input amstex
\documentstyle{amsppt}

\subjclassyear{2000}

\NoBlackBoxes

\def\xspace{\futurelet\dalsitoken\mojxspace}
\def\mojxspace{%
  \ifx\dalsitoken\bgroup\else
  \ifx\dalsitoken\egroup\else
  \ifx\dalsitoken\/\else
  \ifx\dalsitoken\ \else
  \ifx\dalsitoken~\else
  \ifx\dalsitoken.\else
  \ifx\dalsitoken!\else
  \ifx\dalsitoken,\else
  \ifx\dalsitoken:\else
  \ifx\dalsitoken;\else
  \ifx\dalsitoken?\else
  \ifx\dalsitoken/\else
  \ifx\dalsitoken'\else
  \ifx\dalsitoken)\else
  \ifx\dalsitoken-\else
  \ifx\dalsitoken\space\else
   \space\ignorespaces
   \fi\fi\fi\fi\fi\fi\fi\fi\fi\fi\fi\fi\fi\fi\fi\fi}
\def\ensuremath#1{\ifmmode#1\else$#1$\fi}
\def\emph#1{{\it #1}}
\def\endd{\qed\enddemo}
\def\citeAHS{\cite{1}\xspace}
\def\citeARHIT{\cite{2}\xspace}

\def\citeCINC{\cite{3,Example 4.9}\xspace}
\def\citeCINB{\cite{3,Proposition 3.1}\xspace}
\def\citeCINA{\cite{3,Proposition 3.5}\xspace}

\def\citeENG{\cite{5,Theorem 2.3.26}\xspace}
\def\citeGEOR{\cite{6}\xspace}
\def\citeHER{\cite{7}\xspace}
\def\citeHERHUS{\cite{8}\xspace}

\def\citeKANN{\cite{9,Remark 2.4.4(5)}\xspace}
\def\citeNYI{\cite{10}\xspace}
\def\citeZHOU{\cite{11}\xspace}
\def\refSAISCOR{2.1\xspace}
\def\refPROPPRIMF{2.2\xspace}
\def\refSEKVENAUZ{3.1\xspace}
\def\refSWJEBETA{3.2\xspace}
\def\refGENALASPSRA{3.3\xspace}
\def\refPODPRPSR{3.4\xspace}

\def\refPROPSALJ{3.5\xspace}
\def\refTHMTOPJ{3.6\xspace}
\def\refPROPOSALP{4.1\xspace}
\def\refCOLINT{4.2\xspace}
\def\refCOL{4.3\xspace}
\def\refPROPMALRET{4.4\xspace}
\def\refNOVAVETA{4.5\xspace}
\def\refMALPJESALP{4.6\xspace}
\def\refMALJESAL{4.7\xspace}
\def\refOBSMAL{4.8\xspace}
\def\refREMTIGH{4.9\xspace}

\def\Kat#1{\ensuremath{\bold{#1}}\xspace}
\def\Top{\Kat{Top}}
\def\FG{\Kat{FG}}
\def\SSb#1{\ensuremath{\roman S #1}\xspace}
\def\SSp#1{\SSb{\Kat{#1}}}
\def\CHb#1{\ensuremath{\roman{CH(}#1\roman{)}}\xspace}
\def\CH#1{\CHb{\Kat{#1}}}

\def\Genb#1{\ensuremath{\Kat{Gen}(#1)}\xspace}
\def\Genal{\Genb{\alpha}}
\def\PsRa#1{\ensuremath{\Kat{Psrad}(#1)}\xspace}
\def\PsRad{\Kat{PsRad}}
\def\SPsRa#1{\ensuremath{\Kat{SPsrad}(#1)}\xspace}

\def\On{\ensuremath{\roman{On}}}
\def\Cn{\ensuremath{\roman{Cn}}}
\def\RCn{\ensuremath{\roman{RCn}}}


\def\Calp{\ensuremath{C(\alpha)}}
\def\Malp{\ensuremath{M(\alpha)}}
\def\Balp{\ensuremath{B(\alpha)}}
\def\Cbet{\ensuremath{C(\beta)}}
\def\Mbet{\ensuremath{M(\beta)}}
\def\Salj{\ensuremath{(S^\alpha)_1}}
\def\inv#1{\ensuremath{{#1}^{-1}}}
\def\Invobr#1#2{\ensuremath{\inv{#1}(#2)}}
\def\Obr#1#2{\ensuremath{#1(#2)}}
\def\topo#1{\ensuremath{\Cal{#1}}}
\def\Vloz#1#2#3{\ensuremath{#1\colon #2\hookrightarrow #3}}
\def\Zobr#1#2#3{\ensuremath{#1\colon #2\to #3}}
\def\card#1{\ensuremath{\operatorname{card}#1}}
\def\Ldots#1#2#3{\ensuremath{#1_{#2},\ldots,#1_{#3}}}
\def\eg{e.g.\ }
\def\ie{i.e.\ }
\def\Ie{I.e.\ }
\def\iaoi{if and only if\xspace}
\def\tsp{topological space\xspace}
\def\tsps{topological spaces\xspace}

\def\primfac{prime factor\xspace}

\def\ssp{subspace\xspace}
\def\Sierp{Sierpi\'nski\xspace}
\def\sekven#1{$#1$-sequential\xspace}

\hsize=124mm \vsize=186mm

\topmatter
\title
Subspaces of pseudoradial spaces
\endtitle
\abstract We prove that every topological space ($T_0$-space,
$T_1$-space) can be embedded in a pseudoradial space (in a
pseudoradial $T_0$-space, $T_1$-space). This answers the Problem 3
in \citeARHIT. We describe the smallest coreflective subcategory
\Kat A of \Top such that the hereditary coreflective hull of \Kat
A is the whole category \Top.

\endabstract
\author
Martin Sleziak
\endauthor
\address Department of Algebra and Number Theory, FMFI UK, Mlynsk\'a dolina, 842 48 Bratislava\endaddress
\email \tt sleziak\@fmph.uniba.sk\endemail
\subjclass Primary 54B30; Secondary 18B30 \endsubjclass
\keywords pseudoradial space, coreflective subcategory, coreflective hull, prime space, $\beta$-sequential space
\endkeywords
\endtopmatter

\document

\head 1. Introduction \endhead

Pseudoradial or chain-net spaces were introduced by H.\ Herrlich
in \citeGEOR. In the paper \citeARHIT A.\ V.\
Arhangel'ski\u\i, R.\ Isler and G.\ Tironi asked whether every
topological space is a subspace of a pseudoradial space. The
question was asked again in Nyikos' survey \citeNYI.

In \citeZHOU J.\ Zhou proved that under the assumption
$\frak p=\frak c$  every
countable $T_2$-prime space is a subspace of a pseudoradial
$T_2$-space and, as a consequence, he obtained that every space
of countable tightness embeds in a pseudoradial space.

In this paper we show that every topological space ($T_0$-space,
$T_1$-space) can be embedded in a pseudoradial space (in a
pseudoradial $T_0$-space, $T_1$-space). This follows from the fact
that any topological power of the \Sierp doubleton is a
pseudoradial \hbox{($T_0$-)}space.

We also give a characterization of coreflective subcategories $\Kat A$ of $\Top$ for which
every space can be embedded in a space that belongs to $\Kat A$.

\head 2. Preliminaries and notations \endhead

The classes of spaces investigated in this paper are closed under the formation of
topological sums and quotient spaces. In the categorical language, they are coreflective
subcategories of the category \Top of topological spaces. We recall some properties of
coreflective subcategories of \Top which seem to be useful for our investigations (see
\eg \citeHER, \citeAHS).

Let $\Kat A$ be a full and isomorphism-closed subcategory of \Top. Then $\Kat A$ is
\emph{coreflective} \iaoi it is closed under the formation of topological sums and
quotient spaces. If $\Kat B$ is a class of topological spaces (a subcategory of $\Top$)
then by $\CH B$ we denote the \emph{coreflective hull} of $\Kat B$ \ie the smallest
coreflective subcategory of $\Top$ containing $\Kat B$. $\CH B$ consists of all
quotients of sums of spaces which belong to $\Kat B$.

Let $\Kat A$ be a subcategory of $\Top$ and let  $\SSp A$ denote the subcategory of $\Top$
consisting of all subspaces of spaces from $\Kat A$. Then the following result is known (see
\citeKANN or \citeCINB).

\proclaim{Proposition \refSAISCOR} If \Kat A is a coreflective subcategory of $\Top$, then
$\SSp A$ is also a coreflective subcategory of $\Top$. ($\SSp A$ is the coreflective
hereditary hull of $\Kat A$.)
\endproclaim

Given a \tsp $X$ and a point $a\in X$, denote by $X_a$ the space constructed
by making each point, other then $a$, isolated with $a$ retaining its original
neighborhoods. (\Ie a subset $U\subseteq X$ is open in $X_a$ \iaoi $a\notin U$
or there exists an open subset $V$ of $X$ such that $a\in V\subseteq U$.)

We say that a coreflective subcategory \Kat A of \Top is
\emph{nontrivial} if $\FG\subseteq\Kat A$. ($\FG$ denotes the
class of all finitely generated \tsps.) In \citeCINA it is shown
that

\proclaim{Proposition \refPROPPRIMF} If $\Kat A$ is a nontrivial hereditary coreflective
subcategory of $\Top$, then for each $X\in\Kat A$ and each $a\in X$ the \primfac $X_a$ of $X$
at $a$ belongs to $\Kat A$.
\endproclaim

Cardinals are initial ordinals where each ordinal is the
(well-ordered) set of its predecessors. We denote the class of all
ordinal numbers by $\On$, the class of all infinite cardinals by
$\Cn$ and the class of all regular cardinals by $\RCn$.

\emph{Transfinite sequence} is a net defined on an infinite ordinal. In particular, a
transfinite sequence defined on the ordinal $\alpha$ is said to be an
\emph{$\alpha$-sequence.}

A \tsp $X$ is said to be a \emph{prime space} if it contains precisely one accumulation
point.

Finally, let $t(X)$ denote the tightness of $X$ and $\alpha$ be an infinite
cardinal. By $\Genal$ we denote the subcategory of $\Top$ consisting of all
spaces $X$ with $t(X)\leq\alpha$. It is well known that $\Genal$ is a
coreflective subcategory of $\Top$. Moreover, it is the coreflective hull of
the class of all prime spaces $P$ with $\card P\leq\alpha$.

\head 3. Subspaces of pseudoradial spaces \endhead

We start with the definition of pseudoradial and \sekven\beta
space.

A topological space $X$ is said to be \emph{pseudoradial} if, for
any subset $A$ of $X$, $A$ is closed whenever together with any
transfinite sequence it contains all its limits. Let $\beta$ be
an infinite cardinal. A space $X$ is said to be
\emph{\sekven\beta} if, for any subset $A$ of $X$, $A$ is closed
whenever together with any $\alpha$-sequence such that
$\alpha\leq\beta$ it contains all its limits.

Observe, that if $X$ is a \sekven\beta space, then $X$ is
pseudoradial and if $\beta\leq\gamma$ and $X$ is \sekven\beta
then $X$ is \sekven\gamma.

It is useful to characterize \sekven\beta spaces using
\sekven\beta closure. Let $X$ be a topological space and
$A\subseteq X$. The \emph{\sekven\beta closure} of $A$ is the
smallest set $\widetilde A$ such that $A\subseteq \widetilde A$
and $\widetilde A$ is closed with respect to limits of
$\alpha$-sequences for every $\alpha\leq\beta$. Obviously, if $A$
is a subset of $X$, then $\widetilde A \subseteq \overline A$ and
if $A\subseteq B\subseteq X$, then $\widetilde A \subseteq
\widetilde B$. The following characterization of \sekven\beta
spaces is well known and easy to see.

\proclaim{Proposition \refSEKVENAUZ}
A \tsp $X$ is \sekven\beta \iaoi for any subset $A$ of $X$
$\widetilde A = \overline A$.
\endproclaim

We denote by $\PsRad$ the (full) subcategory of $\Top$ consisting of all pseudoradial spaces,
by $\PsRa\beta$ the subcategory consisting of all \sekven\beta spaces and by $\SPsRa\beta$ the
subcategory of all subspaces of \sekven\beta spaces. It is well known that $\PsRad$ and
$\PsRa\beta$ are coreflective subcategories of $\Top$ and, consequently, $\SPsRa\beta$ is also
coreflective in $\Top$.

Denote by $\Calp$ the \tsp on $\alpha\cup\{\alpha\}$ such that a
subset $U\subseteq \alpha\cup\{\alpha\}$ is open in $\Calp$ \iaoi
$U\subseteq\alpha$ or $\card{(\Calp\setminus U)}<\alpha$. It is
known that $\PsRa\beta=\CHb{\{\Calp; \alpha\leq\beta;
\alpha\in\RCn\}}$ and $\PsRad=\CHb{\{\Calp; \alpha\in\RCn\}}$.

Next we want to prove that for any infinite cardinal $\alpha$
$\Genal \subseteq \SPsRa{2^\alpha}$. As a consequence we obtain
that every topological space is a subspace of a pseudoradial
space.

Denote by $S$ the \Sierp space, \ie the space defined on the set
$\{0,1\}$ with the topology consisting of the empty set, the set
$\{0\}$ and the whole space.

\proclaim{Proposition \refSWJEBETA}
If $\beta$ is an infinite cardinal, then the topological power $S^\beta$ of the space $S$ is a
\sekven{\beta} space.
\endproclaim

\demo{Proof}
Let $\gamma$ be the smallest cardinal such that $S^\gamma$ is not \sekven\beta. We want to
show that $\gamma>\beta$. Since for any cardinal $\alpha\leq\omega_0$ $S^\alpha$ is a
sequential space (it is first-countable), $\gamma>\omega_0$. Assume that $\gamma\leq\beta$.
According to Proposition \refSEKVENAUZ there exists a subset $U$ of $S^\gamma$ with
$\overline U \setminus \widetilde U \neq \emptyset$ (by $\widetilde U$ we denote the
\sekven\beta closure of $U$).

Let $t\in \overline U \setminus \widetilde U$,
$A=\{\eta\in\gamma; t(\eta)=0\}$ and  $\varkappa=\card A$.
Clearly, $A\neq\emptyset$ and $\varkappa\leq\gamma$. Consider the
subspace $K=\{s\in S^\gamma;$ for each $\eta\in\gamma\setminus A\
s(\eta)=1\}$. The space $K$ is a closed subspace of $S^\gamma$,
$t\in K$ and, obviously, $K$ is homeomorphic to the space
$S^\varkappa$.

Let us define a map $\Zobr g{S^\gamma}K$ by
$$
g(f)(x)=
\cases
f(x), &\text{if } x\in A\\
1, &\text{if } x\notin A
\endcases
$$
The map $g$ is continuous, hence $t=g(t)\in \overline{\Obr gU}$. Clearly, for
every $s\in S^\gamma$ $g(s)\in\widetilde{\{s\}}$ (the constant sequence
$(s)_{n<\omega_0}$ converges to $g(s)$) and therefore $\widetilde{\Obr gU}
\subseteq \widetilde U$.

If $\varkappa<\gamma$, then $K$ is \sekven\beta and therefore $\overline{\Obr
gU}=\widetilde{\Obr gU}$. This implies that $t\in\widetilde U$. Thus, we
obtain that $\varkappa = \gamma$.

In this case there exists a homeomorphism $\Zobr fK{S^\gamma}$ such that
$f(t)=t_0$ where $t_0(\eta)=0$ for each $\eta\in\gamma$. Without loss of
generality we can suppose that $K=S^\gamma$ and $t=t_0$ (and, obviously, $g$
is the identity map). For each $\xi\in \gamma$ let $f_\xi$ denote the element
of $S^\gamma$ given by
$$
f_\xi(x)=
\cases
0,&\text{for } x<\xi,\\
1,&\text{for } x\geq\xi.
\endcases
$$
It is easy to see that the $\gamma$-sequence
$(f_\xi)_{\xi\in\gamma}$ converges to $t_0$ in $S^\gamma$. Since
$t_0\in \overline U$ and $\overline{\{t_0\}}=S^\gamma$, we obtain
that $\overline U = S^\gamma$ and therefore $f_\xi\in\overline U$
for each $\xi\in\gamma$. Put $A_\xi=\{\eta\in\gamma;
f_\xi(\eta)=0\}=\{\eta\in\gamma;\eta<\xi\}$. Then for each
$\xi\in\gamma$ $\card{A_\xi} < \gamma$ and according to the
preceding part of proof (the case $\varkappa < \gamma$) $f_\xi\in
\widetilde U$. Hence, $t_0 \in \widetilde U$ contradicting our
assumption. Thus, $\gamma > \beta$ and $S^\beta$ is \sekven\beta.
\endd

\proclaim{Theorem \refGENALASPSRA}
$\Genal\subseteq\SPsRa{2^\alpha}$ for every infinite cardinal $\alpha$.
\endproclaim

\demo{Proof}
Since $\Genal$ is the coreflective hull of the class of all prime
spaces $P$ with $\card P \leq \alpha$ and $\SPsRa{2^\alpha}$ is
coreflective, it suffices to prove that every prime space $P$ with
$\card P\leq \alpha$ belongs to $\SPsRa{2^\alpha}$.

Let $P$ be a prime space and $\card P \leq \alpha$. Then $P$ is a $T_0$-space
and the weight of $P$ $w(P) = \beta \leq 2^\alpha$. It is well known (see \eg
\citeENG) that $P$ is embeddable in $S^\beta$. According to Proposition
\refSWJEBETA $S^\beta$ is \sekven\beta and therefore it is also
$2^\alpha$-sequential. Hence $P$ belongs to $\SPsRa{2^\alpha}$.
\endd

As a consequence of the preceding theorem we obtain:

\proclaim{Theorem \refPODPRPSR}
Any \tsp is a subspace of a pseudoradial space. Moreover, every $T_0$-space is a subspace of
a pseudoradial $T_0$-space.
\endproclaim

\demo{Proof}
The first part is an easy consequence of Theorem \refGENALASPSRA.
The second part follows from the
fact that every $T_0$-space is subspace of some $S^\alpha$ (\citeENG)
and from Proposition \refSWJEBETA.
\endd

We next show that this result holds also in the class of all $T_1$-spaces.

Let us recall that the \emph{cofinite topology} on an underlying set $X$ is
the coarsest $T_1$ topology on this set. Closed sets in the cofinite topology
are precisely finite sets and the whole set $X$.

For any cardinal number $\alpha$, let $\Salj$ be the \tsp on the
set $\{0,1\}^\alpha$ with the topology which is the join of the
product topology $S^\alpha$ and the cofinite topology on the set
$\{0,1\}^\alpha$. If $\alpha$ is finite, then $\Salj$ is discrete
space.

\proclaim{Proposition \refPROPSALJ} Let $\alpha$ be an infinite
cardinal. The \tsp $\Salj$ is \sekven\alpha.
\endproclaim

\demo{Proof} The collection
$$\topo B_1=\{U_M; M\subseteq \alpha, M\text{ is finite}\},\text{ where }
U_M=\{f\in \{0,1\}^\alpha; f(m)=0\text{ for each }m\in M\}$$ is
the canonical base for the product topology $S^\alpha$. Clearly
$$\topo B=\{U_M\setminus F; M\subseteq\alpha, M\text{ is finite,}
F\subseteq \{0,1\}^\alpha, F\text{ is finite}\}$$
is a base for the topology of the space $\Salj$.

We have to show that if $t\in\overline U\setminus U$ then
$t\in\widetilde U$. (By $\widetilde U$ we denote the \sekven\alpha
closure of $U$ in $\Salj$.) Let us put
$$A_t=\{\eta\in\alpha; t(\eta)=0\}.$$

Assume that, on the contrary, there exist some $t\in
\{0,1\}^\alpha$ and $U\subseteq \{0,1\}^\alpha$ such that
$t\in\overline U \setminus U$ and $t\notin\widetilde U$. Let
$\beta$ be the smallest cardinal number such that $\beta=\card
A_t$ for some $t$ and $U$ satisfying $t\in\overline U \setminus U$
and $t\notin\widetilde U$.

First let $\beta$ be finite, \ie let $A_t$ be a finite subset of $\alpha$.
Then $U_{A_t}$ is a neighborhood of $t$, thus there exists $f_1\in U\cap
U_{A_t}$. Since $U_{A_t}\setminus\{f_1\}$ is a neighborhood of $t$, there is
$f_2 \in U\cap(U_{A_t}\setminus\{f_1\})$. In a similar way we can find for
every $n<\omega$, $n\geq 2$, an $f_n\in U\cap (U_{A_t}\setminus\{\Ldots
f1{n-1}\})$. We claim that $f_n$ converges to $t$.

Every basic neighborhood of $t$ has the form $U_B\setminus F$,
where $F\subseteq S^\alpha$ and $B\subseteq A_t$ are finite
subsets. $U_B$ contains all terms of the sequence
$(f_n)_{n<\omega}$ and by omitting the finite subset $F$ we omit
only finitely many of them, since this sequence is one-to-one.

Thus $\beta$ is not finite and $\omega \leq \beta = \card A_t \leq \alpha$.
Let us arrange all members of $A_t$ into a one-to-one $\beta$-sequence. Hence,
$A_t=\{a_\xi; \xi<\beta\}$. Let us define a

function $\Zobr{f_\gamma}{\alpha}S$ by
$$
f_\gamma(x)=
\cases
0, &\text{if }x=a_\xi\text{ for some }\xi<\gamma\\
1, &\text{otherwise,}
\endcases
$$
for every $\gamma<\beta$.

If $U_B \setminus F$ is a basic neighborhood of $f_\gamma$, then
$(U_B \setminus F) \cup \{t\}$ is a neighborhood of $t$. Hence
$f_\gamma\in\overline U$. Since the cardinality of the set
$A_\gamma=\{a_\xi, \xi<\gamma\}= \{\eta\in\beta;
f_\gamma(\beta)=0\}$ is less then $\beta$ and $f_\gamma \in
\overline U$, we get $f_\gamma \in \widetilde U$.

It only remains to show that the sequence $f_\gamma$ converges to
$t$. Any basic neighborhood of $t$ has the form $U_B\setminus F$,
where $B\subseteq A_t$, $B$ and $F$ are finite. Let
$\delta_1=\sup\{\xi: a_\xi\in B\}$ and $\delta_2=\sup\{\xi:
f_\xi\in F\}$. Since $F$ and $B$ are finite,
$\delta_1,\delta_2<\beta$. Let
$\delta=\max\{\delta_1,\delta_2\}$. Then for each $\gamma>\delta$
$f_\gamma \in U_B \setminus F$.

Thus $t\in\widetilde U$, a contradiction.
\enddemo

\proclaim{Theorem \refTHMTOPJ} Every $T_1$-space is a \ssp of a
pseudoradial $T_1$-space.
\endproclaim

\demo{Proof} Let $X$ be a $T_1$-space. Then there exists an
embedding $\Vloz eX{S^\alpha}$ of $X$ into some topological power
$S^\alpha$ of $S$. Since $X$ is $T_1$, $\Vloz eX{\Salj}$ is an
embedding as well. $\Salj$ is a $T_1$-space and it is pseudoradial
by Proposition \refPROPSALJ.
\enddemo

\head 4. Coreflective subcategories with $\SSp A=\Top$ \endhead

In \citeHERHUS H.\ Herrlich and M.\ Hu\v{s}ek suggested to investigate the
coreflective subcategories of \Top for which the coreflective hereditary
kernel is the category \FG and the coreflective hereditary hull is the whole
category $\Top$. According to Theorem \refPODPRPSR and \citeCINC the category
$\PsRad$ is an example of such category. Let $\Cal S$ denote the collection of
all such subcategories of $\Top$. We next show that the intersection of any
nonempty family of elements of $\Cal S$ belongs to $\Cal S$ and $\Cal S$ has
the smallest element.

Recall (see Proposition \refSAISCOR) that if $\Kat A$ is a coreflective subcategory of $\Top$,
then the coreflective hereditary hull of $\Kat A$ is $\SSp A$. We first give a
characterization of coreflective subcategories of $\Top$ for which $\SSp A=\Top$.

\proclaim{Theorem \refPROPOSALP} Let $\Kat A$ be a coreflective subcategory of $\Top$. Then
$\SSp A=\Top$ \iaoi $S^\alpha \in \Kat A$ for every infinite cardinal $\alpha$.
\endproclaim

\demo{Proof} Let $\Kat A$ be a coreflective subcategory of $\Top$ for which
$\SSp A=\Top$ and $\alpha$ be any infinite cardinal. There exists a space
$X\in\Kat A$ such that $S^\alpha$ is a subspace of $X$. For each $a\in\alpha$
let $\Zobr{p_a}{S^\alpha}S$ denote the $a$-th projection of topological power
$S^\alpha$ onto $S$. The set $\inv{(p_a)}(0)$ is open in $S^\alpha$ so that
there exists an open subset $U_a$ in $X$ such that $U_a\cap S^\alpha =
\inv{(p_a)}(0)$. The map $\Zobr{f_a}XS$ given by $f_a(x)=0$ for each $x\in
U_a$ and $f_a(x)=1$ otherwise is a continuous extension of
$\Zobr{p_a}{S^\alpha}S$. The map $\Zobr fX{S^\alpha}$ with $f_a=p_a\circ f$
for each $a\in\alpha$ is continuous and the restriction $f|_{S^\alpha}$ is the
identity map on $S^\alpha$. Hence $f$ is a retraction and, consequently, $f$
is a quotient map. Thus $S^\alpha \in \Kat A$.

Conversely, if for any cardinal $\alpha$ $S^\alpha$ belongs to $\Kat A$, then any prime space
belongs to $\SSp A$ and since $\SSp A$ is a coreflective subcategory of $\Top$ we obtain that
$\SSp A=\Top$.
\endd

\proclaim{Corollary \refCOLINT} If $\{\Kat A_i, i\in I\}$ is a nonempty collection of
coreflective subcategories of $\Top$ such that for each $i\in I$ $\SSp A_i=\Top$ and $\Kat
A=\bigcap\{\Kat A_i, i\in I\}$, then $\SSp{A}=\Top$.

If, moreover, for each $i\in I$ the coreflective hereditary kernel of $\Kat A_i$ is $\FG$,
then, obviously, the coreflective hereditary kernel of $\Kat A$ is again $\FG$.
\endproclaim

\proclaim{Corollary \refCOL} $\Kat A=\CHb{\{S^\alpha; \alpha\in\Cn\}}$ is the
smallest coreflective subcategory of \Top such that $\SSp A=\Top$. Obviously,
the coreflective hereditary kernel of $\Kat A$ is $\FG$ (since $\FG \subseteq
\Kat A \subseteq \PsRad$ and $\FG$ is the coreflective hereditary kernel of
$\PsRad$).
\endproclaim

Note that Theorem \refPROPOSALP, Corollary \refCOLINT and Corollary \refCOL remain valid
after replacing $\Top$ by $\Top_0$ (the category of $T_0$-spaces).

We next present another class of (in some sense more convenient) generators of
the category $\CHb{\{S^\alpha; \alpha\in\Cn\}}$.

Let $\alpha$ be an infinite cardinal and
$B_\beta=\{\gamma\in\alpha\cup\{\alpha\}; \gamma\geq\beta\}$ for each
$\beta\in\alpha$. Then $\Malp$ is the \tsp on the set $\alpha\cup\{\alpha\}$
with the topology consisting of all $B_\beta$, $\beta$ being a non-limit
ordinal less then $\alpha$ or $\beta=0$. These spaces have the following
useful property:

\proclaim{Proposition \refPROPMALRET} Let $\alpha$ be an infinite cardinal and
$\Malp$ be a \ssp of $X$. Then there exists a retraction $\Zobr fX\Malp$.
\endproclaim

\demo{Proof} For every non-limit ordinal $\beta<\alpha$ denote by
$U_\beta$ the union of all open subsets of $X$ with
$U\cap\Malp=B_\beta$ and put $U_0=X$. Clearly, if $0\leq \beta <
\beta'<\alpha$ then $U_\beta \supsetneqq U_{\beta'}$ and for each
$\beta<\alpha$ $U_\beta \cap \Malp = B_\beta$. Define $\Zobr
fX\Malp$ by
$$f(x)=\sup\{\beta\in\alpha: x\in U_\beta\}.$$
Obviously, $\Invobr f{B_\beta}=U_\beta$ for non-limit ordinal
$\beta$. Thus $f$ is continuous. Moreover we have $f(\beta)=\beta$
for $\beta\in\Malp$ and $f$ is a retraction.
\endd

\proclaim{Theorem \refNOVAVETA} Let $\Kat A$ be a coreflective subcategory of
$\Top$ and $\alpha$ be an infinite cardinal. The following statements are
equivalent: \roster
\item $\PsRa\alpha \subseteq \SSp A$
\item $S^\alpha \in \Kat A$
\item $\Malp\in\Kat A$
\endroster
\endproclaim

\demo{Proof} $(1)\Rightarrow (2)$ By Proposition \refSWJEBETA
$S^\alpha\in\PsRa\alpha$. Hence $S^\alpha \in \SSp A$, \ie $S^\alpha$ is a
\ssp of a space $X\in\Kat A$. Following the proof of Theorem \refPROPOSALP
we can construct a retraction $\Zobr fX{S^\alpha}$, thus $S^\alpha\in\Kat A$.

$(2)\Rightarrow(3)$ Let $S^\alpha\in\Kat A$. The weight of the space $\Malp$ is
$w(\Malp)=\alpha$, therefore $\Malp$ is a \ssp of $S^\alpha$ by \citeENG. Then
by Proposition \refPROPMALRET there exists a retraction $\Zobr
g{S^\alpha}\Malp$ and $\Malp\in\Kat A$.

$(3)\Rightarrow(1)$ Let $\Malp\in\Kat A$. Clearly, $\Mbet$ is a \ssp of $\Malp$ for every
$\beta<\alpha$. ($\Mbet$ is the \ssp on the set $\beta\cup\{\beta\}$.) Thus for every
$\beta\leq\alpha$ we have $\Mbet\in\SSp A$ and $\Cbet=(\Mbet)_\beta \in \SSp A$ (using
Proposition \refPROPPRIMF). Therefore $\PsRa\alpha = \CHb{\{\Cbet; \beta\leq\alpha\}}
\subseteq \SSp A$.
\endd

\proclaim{Corollary \refMALPJESALP} Let $\alpha$ be an infinite cardinal number. Then
$\CHb{\Malp} = \CHb{S^\alpha}$ and this is the smallest coreflective subcategory of $\Top$
such that $\PsRa\alpha\subseteq\SSp A$.
\endproclaim

\proclaim{Corollary \refMALJESAL} $\CHb{\{\Malp; \alpha\in\RCn\}}=\CHb{\{S^\alpha;
\alpha\in\Cn\}}.$
\endproclaim

\proclaim{Corollary \refOBSMAL} Let $\Kat A$ be a coreflective subcategory of
$\Top$. Then $\SSp A=\Top$ \iaoi $\Malp \in \Kat A$ for every regular cardinal
$\alpha$.
\endproclaim

For a \tsp $X$ and $x\in X$, $t(X,x)$ denotes the tightness of the point $x$
in the \tsp $X$.

For any infinite cardinal $\alpha$, let $\Balp$ be the \tsp on the set
$\alpha\cup\{\alpha\}$ with the topology consisting of all sets
$B_\beta=\{\gamma\in\alpha\cup\{\alpha\}; \gamma\geq\beta\}$ where
$\beta<\alpha$.

\proclaim{Proposition \refREMTIGH} Let $\Kat A$ be a coreflective
subcategory of $\Top$. Then $\SSp A=\Top$ \iaoi for every regular
cardinal $\alpha$ $\Kat A$ contains a space $X$ such that there
exists $x\in X$ with $t(X,x) = \alpha$ and for $\alpha=\omega_0$
the \primfac $X_x$ of $X$ at $x$ is, moreover, not finitely
generated.
\endproclaim

\demo{Proof} One direction follows from $t(\Malp,\alpha)=\alpha$
and Theorem \refNOVAVETA.

Now let $t(X,x)=\alpha$ and $X\in\Kat A$. Then there exists $V\subseteq X$ and
$x\in X$ with $\card V=\alpha$, $x\in\overline V$ and $x\notin \overline U$
for any $U\subseteq V$, $\card U<\alpha$. Let $Y$ be the \ssp of $X$ on the
set $V\cup\{x\}$. $Y$ belongs to $\SSp A$ and by Proposition \refPROPPRIMF
$Y_x$ also belongs to $\SSp A$.

Next we want to prove that $\Balp\in\SSp A$.

We claim that the \tsp $Y_x$ is finer than $\Balp$. Indeed, if $x\in U
\subseteq Y$ and $\card{(V\setminus U)}<\alpha$ then
$x\notin\overline{V\setminus U}$ and $U$ is neighborhood of $x$ in $Y$, hence
$U$ is open in $Y_x$. Clearly, the set $\{x\}$ is not open in $Y$.

Since $\card Y=\alpha$, we can assume that $Y$ is a \tsp on the set
$\alpha\cup\{\alpha\}$ and $x=\alpha$. For every $\gamma<\delta\leq\alpha$ let
$S_{\gamma\delta}$ be the \Sierp \tsp on the set $\{\gamma,\delta\}$ where the
set $\{\delta\}$ is open. A subset $U\subseteq\alpha\cup\{\alpha\}$ is open in
$\Balp$ \iaoi it is open in $Y_x$ and $U\cap\{\gamma,\delta\}$ is open in
$S_{\gamma\delta}$ for every $\gamma<\delta\leq\alpha$ (\ie $U$ contains with
$\gamma\in U$ every $\delta>\gamma$). Thus $\Balp$ is a quotient space of
$Y_x\sqcup (\coprod S_{\gamma\delta})$ and $\Balp\in\SSp A$.

Then the prime factor $(\Balp)_\alpha=\Calp$ belongs to $\SSp A$ for every
regular cardinal $\alpha$, hence $\PsRad\subseteq\SSp A$ and  by Theorem
\refPODPRPSR $\SSp A=\Top$.
\endd

\noindent{\bf Acknowledgement.} Author would like to thank
J.~\v{C}in\v{c}ura for several comments that led to improvements
in the presentation of results.

\noindent{\bf Note.} The same result without any separation axiom was proved already in the master thesis E. Murtinov\'a: Podprostory pseudoradi\'aln{\'\i}ch prostoru. I wasn't aware of this work when preparing this paper.

\enddocument